\documentclass[12pt,a4paper]{article}
\usepackage[latin1]{inputenc}
\usepackage{amssymb}
\usepackage{pictex}
\usepackage{epsfig}

\newcommand{\bprf}{{\it Proof.~}}
\newcommand{\eprf}{\hfill $\square$ \bigskip\par}

\newcommand{\bn}{\mathcal{B}_n}
\newcommand{\pitr}{ \mathbb{P}_3}

\newcommand{\piu}{\mathbb{P}_1}

\newcommand{\aq}{\tilde{A_4}}
\newcommand{\sq}{\tilde{S_4}}
\newcommand{\ac}{\tilde{A_5}}
\newcommand{\ci}{ \mathbb{C}}

\newcommand{\cix}{\ci[x_0,x_1,x_2,x_3]}
\begin{document}

\begin{center}
{\Large Pencils of symmetric surfaces in $\pitr$}\\[0.3cm]
A.Sarti\\[0.3cm]
June 7, 2001\\[0.3cm]
Mathematisches Institut der Universit\"at\\
Erlangen, Bismarckstra\ss e 1$\frac{1}{2}$\\
Germany\\
e-mail:sarti@mi.uni-erlangen.de\\
\end{center}
\paragraph{0. Introduction.} 
The rotations leaving invariant a platonic solid form a finite subgroup of SO$(3)$. Since the cube is dual to the octahedron, the icosahedron dual to the dodecahedron, there are three such rotation groups $G$$\subseteq$SO$(3)$, isomorphic with $A_4$, $S_4$, and $A_5$ respectively. By {\it binary polyhedral group} we denote the inverse image $\tilde{G}$$\subseteq$SU$(2)$ of $G$ under the universal covering SU$(2)$$\rightarrow$SO$(3)$. The $2:1$ image of $\tilde{G}\times\tilde{G}$ in SO$(4)$ under SU$(2)$$\times$SU$(2)$$\rightarrow$SO$(4)$ will be called {\it bi-polyhedral group} $G_n$:

\begin{center}
\begin{tabular}{l*{3}{l}}
group:&$G_6$&$G_8$&$G_{12}$\\
\hline
order:&$288$&$1152$&$7200$\\
\end{tabular}
\end{center}

The aim of this note is to study the simplest non-trivial $G_n$-invariant polynomials.\\
As subgroups of SO$(4)$ the three groups $G_n$ admit the trivial invariant $Q(x)=x_0^2+x_1^2+x_2^2+x_3^2$. In section {\bf 2.} we compute the Poincar\'e series of the three groups $G_n$ acting on $\cix$ and find that the first non-trivial invariant is homogeneous of degree $n$. This is the reason for choosing the notation $G_n$. In section {\bf 4.} such a non-trivial $G_n$-invariant $S_n$ is given explicitly. The equation
\begin{eqnarray*}
S_n(x)+\lambda Q(x)^{\frac{n}{2}}=0,~~~\lambda\in\mathbb{P}_1(\mathbb{C})
\end{eqnarray*}
defines a pencil of $G_n$-invariant surfaces of degree $n$ in $\mathbb{P}_3(\mathbb{C})$. We compute the singular surfaces in each pencil as well as their singularities.\\
It turns out that each pencil contains, except for the multiple quadric $Q(x)^{\frac{n}{2}}$$=0$, precisely four singular surfaces. Their singularities are ordinary double points and form one $G_n$-orbit. The parameters $\lambda$ for these surfaces and their numbers of singularities are:

\renewcommand{\arraystretch}{1.3}

\begin{center}

\begin{tabular}{l||*{3}{l}*{1}{l||}*{4}{l}}
$n$ & \multicolumn{4}{c||}{$\lambda$}&\multicolumn{4}{c}{singularities}\\
\hline 
6 & $ -1$ & $-\frac{2}{3}$ & $-\frac{7}{12}$ & $-\frac{1}{4}$ & $12$& $48$& $48$& $12$\\
8 & $-1$& $-\frac{3}{4}$& $-\frac{9}{16}$& $-\frac{5}{9}$ & 24& 72& 144& 96\\
12 & $ -\frac{3}{32}$& $-\frac{22}{243}$ & $-\frac{2}{25}$& 0&300& 600& 360& 60\\
\end{tabular}
\end{center} 

\renewcommand{\arraystretch}{1.0}
The degrees of these surfaces and their numbers of singularities are quite large. So most of the computations in this note are done by MAPLE. The results seem to be interesting for the following reasons:
\begin{enumerate}
\item S. Mukai \cite{mukai} identified the quotient spaces $\pitr/G_n$ with the Satake compactification ($n=6,8$), resp. with a certain modification of it ($n=12$), of the moduli spaces of abelian surfaces admitting a certain polarization. Our surfaces descend to these quotients and they should therefore be related to modular forms.
\item The reflection group of the regular four-dimensional 24-cell, resp. 600-cell (cf. \cite{cox} chapter VIII p. 145) contains the group $G_6$, resp. $G_{12}$. The invariants of $G_6$, resp. $G_{12}$ in fact are invariant under this bigger group. The existence of these invariants for the reflection groups is known (cf. \cite{cox0}, \cite{racah}), explicit equations however seem not to have been computed before.
\item V. Goryunov at Europroj `96 observed that $G_{12}$ admits an invariant surface of degree $12$ with 600 nodes. Our explicit computations confirm his announcement. This gives the lower bound $\mu(12)\geq 600$, where $\mu(d)$ denotes the maximal number of nodes on a projective surface of degree $d$. Upper bounds for $\mu(d)$ are given in \cite{var}, \cite{miyaoka}. In particular $\mu(12)\leq 645$ by Miyaoka's bound \cite{miyaoka}. By a result of Krei\ss~\cite{kreiss} it was known that $\mu(12)\geq 576$. Our result explicitly improves this bound to show   
\begin{eqnarray*}
600\leq\mu(12)\leq 645.
\end{eqnarray*}
\end{enumerate}

\noindent
We exhibit in section {\bf 12.} a  computer picture of this surface.\\
\\
{\bf Acknowledgements:} I would like to thank Prof. W. Barth for suggesting me how to find surfaces with many nodes and for many very useful discussions. I also would like to thank Prof. D. van Straten for letting me know about the talk of  V. Goryunov at Europroj `96 and S. Endra\ss~for helping me with the program SURF to draw the computer picture.
\paragraph{1. Rotation groups.}
We start by considering the rotation groups of the platonic solids. These are tetrahedron, octahedron, cube, icosahedron and dodecahedron, where octahedron and cube, icosahedron and dodecahedron are reciprocal (cf. \cite{cox} p. 127), hence they have the same rotation group. We consider tetrahedron, octahedron and icosahedron in the same position as in \cite{cox} p. 52, where for the icosahedron we interchange the axes $x$ and $y$. The rotation groups have representations as matrices of SO$(3)$ and are called {\it polyhedral groups} (cf. \cite{cox} p. 46). We denote them by $T$, $O$ and $I$. These are isomorphic to the permutation groups $A_4$, $S_4$ and $A_5$ respectively ( \cite{cox} p. 46-50, \cite{klein} p. 14-19). We denote by $V$ the Klein four group contained in each group $T,~O,~I$. Consider now the classical surjective $2:1$ maps of \cite{tits} p. 77-78 
\begin{eqnarray*}\label{rho}
\rho:{\rm SU}(2)\rightarrow {\rm SO}(3)~~~\hbox{and}~~~
\sigma:{\rm SU}(2)\times {\rm SU}(2)\rightarrow {\rm SO}(4),
\end{eqnarray*}

\noindent
they transform the polyhedral groups into the {\it binary polyhedral groups} and into subgroups of SO$(4)$ respectively. We call the latter {\it bi-polyhedral groups}. Denote by $\tilde{V}:=\rho^{-1}(V)$, $\aq:=\rho^{-1}(T)$, $\sq:=\rho^{-1}(O)$ and $\ac:=\rho^{-1}(I)$ the binary groups and by $\mathcal{H}:=\sigma(\tilde{V} \times\tilde{V})$, $G_6:=\sigma(\aq \times\aq)$, $G_8:=\sigma(\sq \times\sq)$, $G_{12}:=\sigma(\ac \times\ac)$ the bi-polyhedral groups (these last notations will be clear succesively). Before giving the generators of the groups, we specify some matrices that we will use several times in the sequel. Let $\tau:=\frac{1}{2}(1+\sqrt{5})$.\\
\begin{itemize}
\item In SO$(3)$\\
\begin{eqnarray*}
\begin{array}{ll}
A_1:=\left( \begin{array} {ccc}
1& 0& 0 \\
0 & -1& 0  \\
0 & 0& -1  
\end{array} \right),&
A_2:=\left( \begin{array} {ccc}
-1& 0& 0 \\
0& 1& 0 \\
0 & 0& -1
\end{array} \right),\\
&\\
A_3:=\left( \begin{array} {ccc}
-1& 0& 0 \\
0& -1& 0 \\
0 & 0& 1
\end{array} \right),&
R_3:=\left( \begin{array} {ccc}
0& -1& 0 \\
0 & 0& -1  \\
1 & 0& 0  
\end{array} \right),\\
&\\
R_{4}:=\left( \begin{array} {ccc}
1& 0& 0 \\
0 & 0& -1  \\
0 & 1& 0  
\end{array} \right),&
R_{5}:=\frac{1}{2}\left( \begin{array} {ccc}
\tau-1& -\tau& 1 \\
\tau & 1& \tau-1  \\
-1 & \tau-1& \tau  
\end{array} \right),
\end{array}
\end{eqnarray*}

$A_i^2=R_j^j=id,~i=1,2,3;~j=3,4,5.$

\item In SU$(2)$
the quaternions
\begin{eqnarray*}
\begin{array}{lll}
q_1:=\left( \begin{array} {cc}
i& 0\\
0& -i
\end{array} \right),&
q_2:=\left( \begin{array} {cc}
0& 1\\
-1& 0
\end{array} \right),&
q_3:=\left( \begin{array} {cc}
0& i\\
i& 0
\end{array} \right),
\end{array}
\end{eqnarray*}
$q_i^2=-id,~i=1,2,3$

and the matrices
\begin{eqnarray*}
\begin{array}{ll}
p_3:=\frac{1}{2}\left( \begin{array} {cc}
1+i& -1+i \\
1+i& 1-i
\end{array} \right),&
p_4:=\frac{1}{\sqrt{2}}\left( \begin{array} {cc}
1+i& 0 \\
0& 1-i
\end{array} \right),\\
&\\
p_5:=\frac{1}{2}\left( \begin{array} {cc}
\tau& \tau-1+i \\
1-\tau+i& \tau
\end{array} \right),& p_j^j=-id,~j=3,4,5.\\
\end{array}
\end{eqnarray*}
\item In SO$(4)$\\
\begin{eqnarray*}
\begin{array}{ll}
\sigma_1:=\sigma(q_{1},id)=
\left( \begin{array} {cccc}
0 & -1& 0 & 0 \\
1 & 0& 0 & 0 \\
0 & 0& 0 & -1 \\
0 & 0& 1 & 0 
\end{array}\right),&
\sigma_2:=\sigma(q_{2},id)=
\left( \begin{array} {cccc}
0 & 0& -1 & 0 \\
0 & 0& 0 & 1 \\
1 & 0& 0 & 0 \\
0 & -1& 0 & 0 
\end{array} \right),
\end{array}
\end{eqnarray*}

\begin{eqnarray*}
\begin{array}{ll}
\sigma_3:=\sigma(id,q_{1})=
\left( \begin{array} {cccc}
0 & 1& 0 & 0 \\
-1 & 0& 0 & 0 \\
0 & 0& 0 & -1 \\
0 & 0& 1 & 0 
\end{array} \right),&
\sigma_4:=\sigma(id,q_{2})=
\left( \begin{array} {cccc}
0 & 0& 1 & 0 \\
0 & 0& 0 & 1 \\
-1 & 0& 0 & 0 \\
0 & -1& 0 & 0 
\end{array} \right),
\end{array}
\end{eqnarray*}
\begin{eqnarray*}
\begin{array}{l}
\pi_3:=\sigma(p_3,id)=\frac{1}{2}
\left( \begin{array} {cccc}
1& -1& 1& -1 \\
1 & 1& -1 & -1 \\
-1 & 1& 1 & -1 \\
1 & 1& 1 & 1 
\end{array} \right),\\
\\
\pi_3':=\sigma(id,p_3)=\frac{1}{2}
\left( \begin{array} {cccc}
1 & 1& -1 & 1 \\
-1 & 1& -1 & -1 \\
1 & 1& 1 & -1 \\
-1 & 1& 1 & 1 
\end{array} \right),
\end{array}
\end{eqnarray*}

\begin{eqnarray*}
\begin{array}{l}
\pi_4:=\sigma(p_4,id)=\frac{1}{\sqrt{2}}
\left( \begin{array} {cccc}
1& -1& 0& 0 \\
1 & 1& 0 & 0 \\
0 & 0& 1 & -1 \\
0 & 0& 1 & 1 
\end{array} \right),\\
\\
\pi_4':=\sigma(id, p_4)=\frac{1}{\sqrt{2}}
\left( \begin{array} {cccc}
1 & 1& 0 & 0 \\
-1 & 1& 0& 0 \\
0 & 0& 1 & -1 \\
0 & 0& 1 & 1 
\end{array} \right),
\end{array}
\end{eqnarray*}
\begin{eqnarray*}
\begin{array}{l}
\pi_5:=\sigma(p_5,id)=\frac{1}{2}
\left( \begin{array} {cccc}
\tau & 0& 1-\tau & -1 \\
0 & \tau & -1 & \tau-1 \\
\tau-1 & 1& \tau & 0 \\
1 & 1-\tau& 0 & \tau 
\end{array} \right),\\
\\
\pi_5':=\sigma(id,p_5)=\frac{1}{2}
\left( \begin{array} {cccc}
\tau & 0& \tau-1 & 1 \\
0 & \tau& -1 & \tau-1 \\
1-\tau & 1& \tau & 0 \\
-1 & 1-\tau & 0 & \tau 
\end{array} \right),
\end{array}
\end{eqnarray*}
$\sigma_i^2=\pi_j^j=\pi_j'^j=-id,~i=1,2,3,4;~j=3,4,5.$
\item In O$(4)$\\
\begin{eqnarray*}
\begin{array}{ll}
C=
\left( \begin{array} {cccc}
1 & 0& 0 & 0 \\
0 & -1& 0 & 0 \\
0 & 0& -1& 0 \\
0 & 0& 0 & -1 
\end{array} \right),&
C'=
\left( \begin{array} {cccc}
1 & 0& 0 & 0 \\
0 & 1& 0 & 0 \\
0 & 0& 0& 1 \\
0 & 0& 1 & 0 
\end{array} \right),
\end{array}
\end{eqnarray*}
$C^2=C'^2=id$.
\end{itemize}
 With this notation the group $V$ is generated by $A_1$, $A_2$; $T$ is generated by $A_1$, $A_2$ and $R_3$; $O$ is generated by $A_2$, $R_3$, $R_4$ and $I$ by $A_1$, $A_2$ and $R_5$. The binary groups are generated by the pre-images in SU$(2)$ of the previous generators, hence $\tilde{V}$ is generated by $q_1$, $q_2$; the group  $\aq$ is generated by $q_1$, $q_2$, $p_3$; $\sq$ by $q_2$, $p_3$, $p_4$  and $\ac$ by $q_1$, $q_2$, $p_5$ . Taking the images via $\sigma$ of these matrices we get the generators of the bi-polyhedral groups. We have $\sigma_i$, $i=1,2,3,4$ for $\mathcal{H}$; $\sigma_i$, $i=1,2,3,4$, $\pi_3$, $\pi_3'$ for $G_6$; $\sigma_i$, $i=2,4$, $\pi_j$, $\pi_j'$, $j=3,4$ for $G_8$ and $\sigma_i$, $i=1,2,3,4$, $\pi_5$, $\pi_5'$ for $G_{12}$. These groups have order $\frac{1}{2}|\tilde{G}|^2$, where by $\tilde{G}$ we denote one of the binary groups in SU$(2)$, hence $|G_6|=288$, $|G_8|=1152$, $|G_{12}|=7200$.\\
\\
{\bf (1.1)} {\it Eigenvalues and eigenvectors of the matrices $\sigma(p,q)$}. The tensor pro\-duct $\ci^2\otimes\ci^2$ can be identified with the space of complex $2\times2$-matrices. The simple tensors $v\otimes w$$\in$$\ci^2\otimes\ci^2$, $v=(x_1,x_2)$, $w=(y_1,y_2)$ then are identified with the zero determinant matrices via:
\begin{eqnarray*}\label{tensor}
\begin{array}{lll}
 \ci^2 \times \ci^2 &\longrightarrow & \ci^4\\
((x_1,x_2),(y_1,y_2)) &\longmapsto & \left(\begin{array}{c}
x_1\\
x_2
\end{array}\right)\cdot \left(\begin{array}{cc}
y_1& y_2
\end{array}\right)=
\left( \begin{array} {cc} 
x_1y_1& x_1y_2\\
x_2y_1& x_2y_2
\end{array} \right).
\end{array}
\end{eqnarray*}
We have $\sigma(p,q)$ $v\otimes w$$=pv\cdot(wq^{-1})^{t}$. Using $q^{-1}=\bar{q}^{t}$, we get $\sigma(p,q)$ $v\otimes w$$= pv\otimes \bar{q}w$. By this formula follows:\\
\noindent
{\it The eigenvalues, resp. the eigenvectors, of $\sigma(p,q)$ are the products, resp. the tensor products, of the eigenvalues, resp. of the eigenvectors, of $p$ and $\bar q$.}
\paragraph{2. Poincar\'e series.}
The groups $G_n$$\subseteq$SO$(4)$ act on $\ci^4$ in a natural way, this induces an action on $\ci[x_0,x_1,x_2,x_3]$. We want to calculate the dimension of $\ci[x_0,x_1,x_2,x_3]^{G_n}_j$, the vector space of $G_n$-invariant homogeneous polynomials of degree $j$. We consider the {\it Poincar\'e series}

\renewcommand\arraystretch{1.3}
\smallskip

\begin{eqnarray*}
p(\cix^{G_n}, t)&:=&\sum_{j=0}^{\infty}{t^j \dim\cix_{j}^{G_n}}\\
&&\\
&=&\frac{1}{|G_n|}\sum_{g \in G_n}{\frac{1}{\det(id-g^{-1}t)}}\\
\end{eqnarray*}

\renewcommand\arraystretch{1.0}

by  {\it Molien's theorem} \cite{Ben}, p. 21.\\
Since det$g=1$, we can replace the denominator of the previous expression by the characteristic polynomial of $g$. Moreover, since elements in the same conjugacy class have the same characteristic polynomial
\begin{eqnarray*}
{\bf (2.1)}~~~p(\cix^{G_n}, t)=\frac{1}{|G_n|}\sum \frac{n_g}{\det(g-id\cdot t)} 
\end{eqnarray*}
where the sum runs over all the conjugacy classes of $G_n$ and $n_g$ denotes the number of elements in the conjugacy class of $g$ in GL$(4,\ci)$ (cf. also \cite{burniside} p. 300, for a similar computation). We calculate now the conjugacy classes in $G_n$. We start with the conjugacy classes of $T$, $O$ and $I$, which are those of the permutation groups $A_4$, $S_4$ and $A_5$. We specify a representative and the size of the corresponding conjugacy class below. For convenience we give the conjugacy classes of $V$ as well. Moreover in the case of $T$ we put the elements $R_3$, and $R_3^2$ in bracket since they have the same eigenvalues and so are conjugate in GL$(3,\ci)$:
\renewcommand{\arraystretch}{1.6}
\begin{eqnarray*}
\begin{array}{ll}
{\rm group}&{\rm conj.~classes}\\
V:&id,1; A_2,3;\\
T:&id,1; A_2,3; \{R_3,4; R_3^2,4\};\\
O:&id,1; A_2,3; R_3,8; R_4,6;  R_3R_4,6;\\
I:&id,1; A_2,15; R_5,12; R_5^2,12; R_3,20.
\end{array}
\end{eqnarray*}
\renewcommand{\arraystretch}{1.0}

\noindent
Let now $G$ denote one of the groups $V,~T,~O,~I$$\subseteq$SO$(3)$, and $\tilde{G}$ one of the binary groups $\tilde{V}$, $\aq$, $\sq$, $\ac$. For the conjugacy class of an element $g$$\in$$G$ we have two conjugacy classes in $\tilde{G}=\rho^{-1}(G)$, namely those of $+\tilde{g}$ and $-\tilde{g}$, the two elements in $\rho^{-1}(g)$. If they are equal then $\tilde{g}$ is a traceless matrix. The only traceless matrices in the binary groups are those in the conjugacy classes of $q_2$ or $p_3p_4$.  We can write the product $p_3p_4=\frac{1}{\sqrt{2}}(q_1+q_3)$, so using the multiplication rules for the quaternions $q_1$, $q_2$, $q_3$ (cf. \cite{tits} p. 77) we find in fact $-q_2$$\in$$[q_2]$ and $-p_3p_4\in[p_3p_4]$.\\
\\
\noindent
We give now the conjugacy classes of the binary groups. In the case of $\aq$ and $\sq$ the elements $p_3$ and $-p_3^2$, $-p_3$ and $p_3^2$, resp. $q_2$ and $p_3p_4$ are conjugate in SU$(2)$, in fact they have the same eigenvalues, so we put them in bracket.  We use the notation  $\pm q$ to indicate the representatives $+q$, $-q$. The number given after such a pair of representatives is the size of the corresponding conjugacy classes.
\renewcommand{\arraystretch}{1.6}
\begin{eqnarray*}
\begin{array}{ll}
{\rm group}&{\rm conj.~classes}\\
\tilde V:& \pm id,1; q_2,6;\\
\aq :&\pm id,1; q_2,6;\{ p_3,4; -p_3^2,4\}; \{-p_3,4; p_3^2,4\};\\
\sq:& \pm id,1; \{q_2,6;  p_3p_4,12\}; \pm p_3,8; \pm p_4,6;\\
\ac:&\pm id,1; q_2,30; \pm p_5,12; \pm p_5^2,12; \pm p_3,20.\\
\end{array}
\end{eqnarray*}
\renewcommand{\arraystretch}{1.0}

\noindent
The conjugacy classes of $\mathcal{H}$, resp. $G_n$, $n=6,8,12$, are the images via $\sigma$ of the conjugacy classes of $\tilde{V}\times\tilde{V}$, $\aq\times\aq$, $\sq\times\sq$, $\ac\times\ac$. Their order is the product of the order of the conjugacy classes in the binary groups except for the conjugacy classes of $q_2$ and $p_3p_4$, where we have to divide by two. We will not to list them now. Our aim is to write the series {\bf (2.1)} for the groups $\mathcal{H}$ and $G_n$, so we need just the conjugacy classes in GL$(4,\ci)$. Here the elements with the same eigenvalues have the same conjugacy class so using {\bf (1.1)} we can find them. We give a representative with the size of the respective conjugacy class below. By $\sigma_{24}$ we denote the product of $\sigma_2$ and $\sigma_4$, the other notations have been introduced in the first paragraph, 
\renewcommand{\arraystretch}{1.6}
\begin{eqnarray*}
\begin{array}{ll}
{\rm group}&{\rm conj.~classes}\\
\mathcal{H}:&\pm id,1; \sigma_2,12;  \sigma_{24},18;\\
G_6:&\pm id,1; \sigma_2,12; \sigma_{24},18; \sigma_2\pi_3',96 ;\pm\pi_3,16; \pm\pi_3\pi_3', 64;\\
G_8:&\pm id,1;\sigma_2,36;\sigma_{24},162;\sigma_2\pi_3',288;\sigma_2\pi_4',216;\pm\pi_3,16; \pm\pi_3\pi_3',64;\\
&\pm\pi_3\pi_4',96;\pm\pi_4,12;\pm\pi_4\pi_4',36;\\
G_{12}:&\pm id,1;\sigma_2,60;\sigma_{24},450;\sigma_2\pi_5',720;\sigma_2\pi_5'^2,720;\sigma_2\pi_5'^2\sigma_4,720;\pm\pi_5,24;\\&\pm\pi_5\pi_5',144
\pm\pi_5\pi_5'^2\sigma_4,480;\pm\pi_5^2,24;\pm\pi_5^2\pi_5',288;\pm\pi_5^2\pi_5'^2,144;\\& \pm\pi_5^2\pi_5'^2\sigma_4, 480;\pm\pi_5^2\sigma_2,40;\pm\pi_5^2\sigma_2\pi_5'^2\sigma_4,400.
\end{array}
\end{eqnarray*}
\renewcommand{\arraystretch}{1.0}

\noindent
Now we can calculate the characteristic polynomial of the previous representatives and we write the sums {\bf (2.1)} for $\mathcal{H}$ and $G_n$. A computation with MAPLE shows:
\renewcommand{\arraystretch}{1.6}
\begin{eqnarray*}
\begin{array}{ll}
{\rm  group}&{\rm Poincar\acute{e}~series}\\
\mathcal{H}:&1+t^2+5t^4+6t^6+15t^8+19t^{10}+35t^{12}+44t^{14}+O(t^{16}),\\
G_6:&1+t^2+t^4+2t^6+3t^8+3t^{10}+7t^{12}+8t^{14}+O(t^{16}),\\
G_8:&1+t^2+t^4+t^6+2t^8+2t^{10}+3t^{12}+3t^{14}+O(t^{16}),\\
G_{12}:&1+t^2+t^4+t^6+t^8+t^{10}+2t^{12}+2t^{14}+O(t^{16}).
\end{array}
\end{eqnarray*}
\renewcommand{\arraystretch}{1.0}

\noindent
These series show that:
\begin{itemize}
\item there are invariant polynomials just in even degree (this is due to the well known fact (cf. e.g. \cite{hudson}, \cite{jessop}) that the $\mathcal{H}$-invariant polynomials are just in even degree and $\mathcal{H}$$\subseteq$$G_n$, $n=6,8,12$).

\item In each degree we have the trivial invariant quadric surface $Q_j(x):=(x_0^2+x_1^2+x_2^2+x_3^2)^{\frac{j}{2}}$. In degree $n=6,8,12$ appears the first not trivial $G_n$-invariant polynomial (this explains the notation). We call it $S_n(x)$ and we calculate an expression of it later.
\end{itemize}

\paragraph{3. Reflection groups.}
The group generated by $G_{12}$ and $C$ together (notation of paragraph {\bf 1.}) is the reflection group of order $14400$ (cf. \cite{Ben} p. 80-81), which is the symmetry group of the regular $600$-cell $\{3,3,5\}$ (cf. \cite{cox} p. 153 and \cite{cox0}). The reflection group of order $1152$ which appears on the table  of \cite{Ben} p. 81 is the reflection group of the $24$-cell $\{3,4,3\}$ (cf. \cite{cox} p. 149 and \cite{cox0}), which we take with vertices the permutations of $(\pm1,\pm1,0,0)$ as in \cite{cox} p. 156. The generators of $G_6$ permute these vertices, hence the group is contained in the reflection group of the $24$-cell. We now show\\
{\bf (3.1)} {\it the group $<G_6,C,C'>$ is the reflection group of $\{3,4,3\}$, where $C$, $C'$ are the matrices given in the first paragraph.}\\
\bprf
We have $C^2=id$ and an easy calculation shows that $C\sigma(p,q)C=\sigma(q,p)$, hence $|<C,G_6>|=$$2\cdot 288=$$576$. Clearly $C'$ is a symmetry transformation of the 24-cell. Assume moreover that there are $\sigma$, $\sigma'$ in $G_6$ with $C'\sigma=C\sigma'$ then $CC'=\sigma'\sigma^{-1}$. It follows that $\sigma'\sigma^{-1}$ has order two, hence it is in $\mathcal{H}$. A long calculation shows that we have no matrices in $\mathcal{H}$ which are equal to $CC'$, so $C'G_6\not=CG_6$, $|<G_6,C,C'>|$$=1152$ and the assertion follows.
\eprf
The situation is quite different in the case of $G_8$. It cannot be contained in any symmetry group of a four-dimensional regular polyhedron since the latter are only the 600-cell, the 24-cell and their reciprocals (cf. \cite{cox} p. 145). The group $G_8$ interchanges in fact the $\{3,4,3\}$ and its reciprocal, which we call $\{3,4,3\}'$ (it is a 24-cell again, \cite{cox} p. 149). We take it in a similar position as \cite{cox} p. 156,  with vertices the permutations of $(\pm \sqrt{2},0,0,0)$ and $(\pm \frac{\sqrt{2}}{2},\pm \frac{\sqrt{2}}{2} , \pm \frac{\sqrt{2}}{2}, \pm \frac{\sqrt{2}}{2})$, i.e. it is the reciprocal with respect to the sphere $x_0^2+x_1^2+x_2^2+x_3^2=\sqrt{2}$ (cf. \cite{cox} p. 126). The rotations of order four $\pi_4$, $\pi_4'$ send in fact the vertices of $\{3,4,3\}$ to the vertices of $\{3,4,3\}'$.\\
{\bf (3.2)} {\it Invariant polynomials of the reflection groups}. It is well known (cf. \cite{cox0}) that the invariant polynomials under $<G_{12},C>$ and $<G_6,C,C'>$ are just in even degree. A basic set for the invariants in sense of \cite{cox0} p. 774 is given by polynomials in degrees $2,12,20,30$, respectively $2,6,8,12$ (cf. \cite{racah}, \cite{cox0}). Clearly the invariants under $G_{12}$ are exactly  those of $<G_{12},C>$. In the case of $G_6$ we have to be more careful, because of the matrix $C'$, which describes an odd permutation of the coordinates. In degree six, however it turns out that the groups $G_6$ and $<G_6,C,C'>$ have the same invariant polynomials.


\paragraph{4. Equation of the invariant polynomials.}
We describe briefly the method to calculate $S_n(x)$. We start with a basis of $\ci[x_0,x_1,x_2,x_3]^{\mathcal{H}}_n$, $n=6,8,12$, which is well known (cf. e.g. \cite{hudson}, \cite{jessop}). Then we write a polynomial as a linear combination of the elements of this basis and we impose it to be invariant under the remaining generators of $G_n\backslash\mathcal{H}$. In degree $6$ and $8$ the computations are not difficult if left to MAPLE. In degree $12$ it is however difficult,  because of the relatively high number of parameters involved.\\
We say that a polynomial is {\it totally symmetric} if it is invariant under each coordinate permutation. In degree $6$ and $8$ the invariant polynomials are
\begin{eqnarray*}
\begin{array}{ll}
S_6(x)=\sum  x_i^6+15\sum x_i^2x_j^2x_k^2,~~S_8(x)=\sum x_i^8+14\sum x_i^4x_j^4+168x_0^2x_1^2x_2^2x_3^2,\\
\end{array}
\end{eqnarray*}
where the sums run over all the indices $i,j,k=0,1,2,3$, and if some indices appear together they are different from each other. These polynomials are totally symmetric. Before giving the  expression of $S_{12}(x)$ we introduce some abbreviations. Let $y_i:=x_i^2$ and
\begin{eqnarray*}
\begin{array}{lll}
S_{51}:=\sum y_i^5y_j;&S_{42}:=\sum y_i^4y_j^2;&S_{33}:=\sum y_i^3y_j^3;\\S_{411}:=\sum y_i^4y_jy_k;&S_{321}:=\sum y_i^3y_j^2y_k;&S_{222}:=\sum y_i^2y_j^2y_k^2;\\S_{3111}:=\sum y_i^3y_jy_ky_h;&S_{2211}:=\sum y_i^2y_j^2y_hy_k,
\end{array}
\end{eqnarray*}
where for the sums we use the same conventions as before. The totally symmetric part of $S_{12}(x)$ is 
\begin{eqnarray*}
f_s:=2S_{51}-6S_{42}-12S_{411}+14S_{33}+9S_{321}+348S_{3111}+30S_{222}-270S_{2211}.
\end{eqnarray*}
In a similar way we say that a polynomial $P$ is {\it anti-symmetric} if it is invariant under each even coordinate permutation and $\gamma \cdot P$=$-P$ for each odd permutation $\gamma$. The anti-symmetric part of $S_{12}(x)$ is 
\begin{eqnarray*}
\begin{array}{c}
33\sqrt{5}f_a
\end{array}
\end{eqnarray*}
with
\begin{eqnarray*}
\begin{array}{lll}\\
f_a&:=&y_0^3(y_1^2y_2-y_1y_2^2+y_2^2y_3-y_2y_3^2+y_3^2y_1-y_3y_1^2)-y_1^3(y_2^2y_3-y_2y_3^2\\
&&+y_3^2y_0-y_3y_0^2+y_0^2y_2-y_0y_2^2)+y_2^3(y_0^2y_1-y_0y_1^2+y_1^2y_3-y_1y_3^2\\
&&+y_3^2y_0-y_3y_0^2)-y_3^3(y_0^2y_1-y_0y_1^2+y_1^2y_2-y_1y_2^2+y_2^2y_0-y_2y_0^2).
\end{array}
\end{eqnarray*}
In conclusion
\begin{eqnarray*}\label{s12x}
\begin{array}{c}
S_{12}(x):=f_s+33\sqrt{5}f_a.
\end{array}
\end{eqnarray*}
Observe that $Q_2(x)$ and $S_n(x)$ are algebraically independent. In fact for the point $p=(i\sqrt{2},1,1,0)$ in $\ci^4$ holds $Q_2(p)=0$ but $S_n(p)\not=0$, $n=6,8,12$. This has in particular an interesting consequence. The polynomial $S_8(x)$ is clearly invariant under $<G_6,C,C'>$ and is not the product $S_6(x)Q_2(x)$, hence it is an equation for the non-trivial invariant polynomial of degree eight in the basic set of invariants under $<G_6,C,C'>$.\\
 The homogeneous invariant polynomials define pencils of symmetric surfaces in $\pitr$:
\begin{eqnarray*}
F_n(\lambda):S_n(x)+\lambda Q_n(x)=0,~~~~~~~~~\lambda\in \mathbb{P}_1
\end{eqnarray*}
The aim of the next paragraphs is to investigate the base locus of the pencils and to find the singular surfaces.

\paragraph{5. Base locus.}Let $\mathcal{G}$ be a group acting on $\pitr$.\\
\\
{\bf(5.1) Definition.} A point $z\in\pitr$ is called a {\it fix point} if there is a $\sigma$ $\in$ $\mathcal{G}$ ($\sigma\not=$ $\pm$ $id$), s.t. $\sigma z=z$. We call
\begin{eqnarray*}
{\rm Fix}(z)={\rm Fix}_{\mathcal{G}}(z):=\{g\in \mathcal{G}~|~gz=z\}\subseteq \mathcal{G}
\end{eqnarray*}
{\it fix group} of $z$ and
\begin{eqnarray*} 
O(z)=O_{\mathcal{G}}(z):=\{ gz~|~ g\in \mathcal{G}\}\subseteq\pitr
\end{eqnarray*}
{\it orbit} of $z$. We have the formula:
\begin{eqnarray*}
|{\rm Fix}(z)|\cdot|O(z)|=|\mathcal{G}|
\end{eqnarray*}
{\bf(5.2) Definition.}
We call a line $L\subseteq\pitr$, a {\it line of fix points} (or {\it fix line}) if there is a $\sigma$ $\in$ $\mathcal{G}$ ($\sigma\not=$ $\pm$ $id$) s.t. $\sigma x=x$ for all $x$$\in$$L$.\\
\\
{\bf (5.3)} As in \cite{mukai} we identify $\pitr\backslash Q_2$ with $\mathbb{P}$GL$(2)$, the projective space of invertible complex $2\times2$-matrices. More explicitly  a point $(x_0:x_1:x_2:x_3)$$\in$$\pitr$ corresponds to a matrix \begin{eqnarray*}
x=\left( \begin{array} {cc} 
x_0+ix_1&x_2+ix_3 \\
-x_2+ix_3&x_0-ix_1
\end{array} \right)\in\mathbb{P}{\rm GL}(2)
\end{eqnarray*}  
The quadric $Q_2=\mathbb{P}_1\times\mathbb{P}_1$ corresponds to the rank one complex $2\times 2$-matrices. We use this identification to show:

{\bf(5.4)} {\it The matrices $\sigma(p,id), \sigma(id,q) \in$ $G_n$ have in $\pitr$ two disjoint lines of fix points each. These  are contained in the quadric $Q_2$ and belong to one ruling, respectively to the other ruling of $Q_2$.}\\
\bprf
By {\bf (1.1)} it follows that the matrices $\sigma(p,id), \sigma(id,q)$ have two eigenvalues with multiplicity two each, the eigenspaces are lines of $\mathbb{P}_3$. These are spanned by points which correspond to matrices of rank one in the identification {\bf (5.3)}, hence they are contained in $Q_2$. Again using {\bf (1.1)} we see that these are lines of one ruling resp. of the other ruling of $Q_2$.
\eprf

\noindent
The base locus of the pencil $F_n(\lambda)$ is the variety 
\begin{eqnarray*}
\{x\in\pitr | S_n(x)=Q_n(x)=0\}
\end{eqnarray*}
Clearly it is not reduced and is invariant under the group action. Put $\bn$$:=Q_2\cap S_n$. We consider now the groups $\sigma(\tilde G,id)$ and $\sigma(id, \tilde G)$, ($\tilde G$$=\aq$, $\sq$, $\ac$, as usual) which modulo $\{\pm id\}$ are isomorphic to the subgroups $T$, $O$ and $I$ $\subseteq$ SO$(3)$. It is a well known fact that under the action of these groups there are orbits of the following lengths,

\begin{center}

\renewcommand{\arraystretch}{1.3}

\begin{eqnarray*}
\begin{array}{c|c|c}
{\rm tetrahedron} &{\rm octahedron} &{\rm icosahedron}\\
\hline
12,~6,~4&24,~12,~8,~6&60,~30,~20,~12\\
\end{array}
\end{eqnarray*}

\end{center}

\renewcommand{\arraystretch}{1.0}

\noindent
Moreover, observe that\\
\begin{enumerate}
\item by {\bf (1.1)} it follows that the group $\sigma(\tilde{G},id)$ acts on the lines of the first ruling $\{v\}\times\piu$ and lets invariant each line of the second ruling $\piu\times\{w\}$. Vice versa $\sigma(id, \tilde{G})$ acts on $\piu\times\{w\}$ and lets invariant each line of $\{v\}\times\piu$.
\item  Denote by $\mathcal{L}_n$, $\mathcal{L}_n'$ the sets of lines in $\{v\}\times\piu$, resp. $\piu\times\{w\}$ of the orbit of length $n$. The matrix $C$ maps lines of $\mathcal{L}_n$ to lines of $\mathcal{L}_n'$ (recall that $C\sigma(p,q)C$$=\sigma(p,q)$).
\end{enumerate}
\noindent
Using these facts we show\\
{\bf (5.5)} (a) {\it the variety $\bn$ is reduced, i.e. does not contain multiple components.}\\
(b) {\it The base locus splits in $2n$ lines, $n$ of each ruling of $Q_2$. In particular these are fix lines for elements in $G_n$.}\\

\medskip

\noindent
{\it Proof of ~}(a).
By Bezout's theorem $\deg(Q_2\cap S_n)=2n$. If $Q_2\cap S_n$ is not reduced then there is a component $V\subseteq Q_2\cap S_n$ s.t. $Q_2$ and $S_n$ meet with multiplicity at least two. This is the case when $V$ is singular on $S_n$, or $S_n$ and $Q_2$ are tangent at $V$. Consider a line $L$ of one of the two rulings, not in $\bn$, and which meets $V$ in at least one point. W.l.o.g. assume $L$ in the ruling $\{v\}$$\times$$\piu$. Let $x\in L\cap V$. We have  mult$_x(L\cdot S_n)\geq 2$. The group $\sigma(id, \tilde G)$ acts on L, so we consider the orbit of $x$ under this group. By the table above, we see that $L$ and $ S_n$ meet at more than $n$ points computed with multiplicity, so $L$$\subseteq$$S_n$. This contradicts the assumption. We have shown that $Q_2\cap S_n$ is reduced.
\eprf
\noindent
{\it Proof of ~}(b).
Take a line $L$ of the first or of the second ruling, $L\not\subseteq\bn$. The curve $\bn$ has bi-degree $(n,n)$ on $Q_2$, so $|L\cap\bn|=n$. W.l.o.g. assume that $L$ is in the ruling $\{v\}\times\piu$. The group $\sigma(id,\tilde G)$ acts on the points of $L$. Let $x\in L\cap\bn$, then by the table on the previous page, the orbit of $x$ under $\sigma(id, \tilde{G})$ must have length $n$. Hence $x$ belongs to a line of $\mathcal{L}_n'$. As we have infinitely many lines like $L$, the lines $\mathcal{L}_n'$ are contained in $\bn$. By 2. the lines in $\mathcal{L}_n$ are contained in $\bn$ too.
\paragraph{6. Singular points.}
We give now some general results on the singular points of the pencils.\\
{\bf (6.1)} {\it Let $p$ be a singular point on a surface of the pencil $F_n(\lambda)$ (not on $Q_n$), then $p$ is not contained in the complex quadric.}\\
\bprf
Let $p$ be a singular point on the surface $F_n(\lambda_0)$:=$S_n(x)+\lambda_0Q_n(x)=0$ and assume that $p$$\in$$Q_2$. Then $p\in S_n$, moreover since $\partial_iQ_n(p)=0$, for $i=0,1,2,3$, we have $\partial_iS_n(p)=0$, $i=0,1,2,3$, too, so $p$ is a singular point on $Q_2\cap S_n$. This consists of $2n$ lines which meet each other at $n^2$ points. Hence $p$ must be an intersection point of two lines. If $S_n$ is singular at $p$ it follows that $S_n$ is singular at all the $n^2$ points of intersection of the lines in the base locus, in fact they form one orbit under the action of $\sigma(\tilde{G},id)$ and $\sigma(id, \tilde{G})$ (notation of paragraph {\bf 5.}). In particular $S_n$ has $n$ singular points on a line $L$ in $Q_2\cap S_n$. A hypersurface $\{\partial_i S_n=0\}$ ($i=0,1,2,3$) has degree $n-1$, therefore it intersects $L$  in $n-1$ points. So $S_n$ has at most $n-1$ singular points on $L$. It follows that $L$ is singular on $S_n$. Hence $S_n$ and $Q_2$ meet at $L$ and so at all the $2n$ lines of $Q_2\cap S_n$ with multiplicity at least two. This is not possible, in fact deg$(Q_2\cap S_n)$$=2n$. This shows that $p$$\not\in$$Q_2$.
\eprf
\noindent
This lemma has many consequences:\\
{\bf (6.2)} 
(a) {\it the general surface in the pencil is smooth,} \\
(b) {\it the singular surfaces, not $Q_n$, have only isolated singularities,}\\ 
(c) {\it the surfaces different from $Q_n$ are irreducible and reduced.}\\

 \bprf
(a) It is a consequence of Bertini's theorem.\\
(b) Assume that $S$$:=\{Q_n(x)+\lambda_0S_n(x)=0\}$ contains a singular curve. This meets $Q_2$ in at least one point $p$, which is singular on $S$. By {\bf (6.1)}, this is not possible.\\
(c) Is like (b).
\eprf

{\bf (6.3)}{\it A singular point on a surface in the pencil $F_n(\lambda)$, $n=6,8,12$, is a fix point under $G_n$ in sense of definition {\bf (5.1)}. Moreover, as vector of $\ci^4$, it is eigenvector of a matrix with eigenvalue $+1$ or $-1$.}\\
\bprf 
It is possible to obtain a first rough estimate of the maximal number of singular points on a surface $S$ of degree $n$ in $\pitr$ in the following way: if $S$ has equation $\{F=0\}$ then a point on $S$ is singular if and only if $\partial_iF(p)=0$ for all $i=0,1,2,3$. These surfaces meet in at most $n(n-1)^2$ points with multiplicity. Since these are singular on $S$, they are counted at least two times in the intersection, so the effective bound is $\frac{n}{2}$$(n-1)^2$. We give these numbers for $n=6,8,12$ in the table below. In the second row we give the length of the orbit of a point under $G_n$, which is not a fix point,

\renewcommand{\arraystretch}{1.3}

\begin{center}

\begin{tabular}{c|*{2}{c|}*{1}{c}}
$n$&6&8&12\\
\hline
$\frac{n}{2}$$(n-1)^2$&75&196&726\\
\hline
orbit &144&576&3600
\end{tabular}

\end{center}

\renewcommand{\arraystretch}{1.0}

\noindent
Clearly such a point cannot be singular. Let now $x$ denote a singular point and consider it as vector in $\ci^4$. Let $\sigma:=\sigma(p,q)$$\in$$G_n$ s.t. $\sigma x=\lambda x$ equivalently
\begin{eqnarray*}
pxq^{-1}=\lambda x.
\end{eqnarray*}
Consider $x$ as matrix of $\mathbb{P}$GL(2) as in the identification {\bf (5.3)} and take the determinant on both sides of the previous equation. We get det$(x)=\lambda^2$det$(x)$. In fact det$(p)=$det$(q)$=1 since they are matrices in SU$(2)$. The equality holds only when det$(x)=0$ or $\lambda^2=1$. If det$(x)$=0 then $x\in Q_n$ and this is not possible by {\bf (6.1)}.
\eprf
By this fact follows\\
{\bf (6.4)} {\it the singular points are contained in  fix lines of $G_n$, $n=6,8,12$.}


\paragraph{7. Lines of fix points.}
Let $L_1$, $L_2$, $L_3$, $L_4$ be the fix lines of the elements $\sigma(p,id)$, resp. $\sigma(id,q)$ $\in G_n$. Let $z_{ij}$$:=$$L_i\cap L_j$, $i=1,2$, $j=3,4$, denote the intersection points. If the matrix $\sigma(p,q)$ has the eigenvalue $1$ or $-1$ then it has at least one line, $L$,  of fix points, which form the following configuration with the fix lines of $\sigma(p,id)$ and of  $\sigma(id,q)$


\vspace{-0.5cm}

\begin{center}

\unitlength1.0cm
\begin{picture}(2,4)
\put(-1,0) {\line(1,0){4.0}~~{$L_2$}}
\put(-0.5,0.2){$z_{24}$} 
\put(-1,2) {\line(1,0){4.0}~~{$L_1$}} 
\put(0,-1) {\line(0,1){4.0}~~{$L_4$}}
\put(2,-1) {\line(0,1){4.0}~~{$L_3$}}
\put(-1,-1){\line(1,1){4.0}}
\put(-0.5,2.2){$z_{14}$}
\put(1.5,0.2){$z_{23}$}
\put(1.5,2.2){$z_{13}$}
\put(-0.9,-0.5){$L$}
\end{picture}

\end{center}

\vspace{1.0cm}

\noindent
With this notation we give the following result that we will use later,\\ 

{\bf (7.1)} {\it assume that $L$ meets the base locus of the pencil $F_n(\lambda)$. W.l.o.g. let these intersection points be $z_{13},~z_{24}$ as in the picture. Then for each surface $S\not=Q_n$ in the pencil we have mult$_{z_{ij}}(L\cdot S)=1$, $(i,j)=(1,3)$, $(2,4)$.}\\
\bprf 
By {\bf (6.1)} the points $z_{13},z_{24}$ are smooth points on each surface (not $Q_n$) in the pencil $F_n(\lambda)$. The lines of the two rulings of $Q_2$ which meet at these points, are lines of the base locus, hence are contained in $S$. The tangent space of $S$ at the $z_i$ is the plane spanned by these two lines. Clearly this plane does not contain $L$, hence $L$ cannot be tangent to $S$ at $z_{ij}$.
\eprf
\bigskip

\noindent
{\bf (7.2)} We give below a representative of each of the conjugacy classes in $G_n$ (under the action of $G_n$ itself), which have eigenspaces of dimension two  in $\ci^4$ (and so fix lines in $\mathbb{P}_3$) with eigenvalue $1$ or $-1$. As usual we write $\pm \sigma$ to indicate the representatives $+\sigma$ and $-\sigma$. Clearly these elements have the same fix lines. Observe that the elements in the conjugacy classes of $\pi_3\pi_3'$, $\pi_3^2\pi_3'^2$; $\pi_5\pi_5'$, $\pi_5^2\pi_5'^2$; $\pi_3\pi_3'^2$, $-\pi_3^2\pi_3'$ have two by two the same fix lines, so we consider them together. The elements in the conjugacy classes of $\sigma_{24}$ and $\pi_4\pi_4'$ have the same fix lines too. This follows from the fact that we can write each element in $[\sigma_{24}]$ as the square of an element in $[\pi_4\pi_4']$. Moreover since $C\sigma(p,q)$$C$$=\sigma(q,p)$, the elements in the conjugacy classes of $\sigma_2\pi_3'\pi_4'$ and of $\pi_3\pi_4\sigma_4$ have the same fix lines. Next to each representative we write the number of distinct fix lines in the conjugacy class.
 
\renewcommand{\arraystretch}{1.6}

\begin{eqnarray*}
\begin{array}{ll}
{\rm group}&{\rm fix~lines}\\
G_6:&\sigma_{24}, 18;\pm\pi_3\pi_3',\pm\pi_3^2\pi_3'^2, 16; \pm\pi_3^2\pi_3',\pm\pi_3\pi_3'^2, 16;\\
G_8:&\sigma_{24},\pm\pi_4\pi_4', 18; \pm\pi_3\pi_3', 32; \pi_3\pi_4\pi_3'\pi_4', 72;\pi_3\pi_4\sigma_4, \sigma_2\pi_3'\pi_4', 36;\\
G_{12}:&\sigma_{24}, 450; \pm \pi_3\pi_3',200; \pm\pi_5\pi_5', \pm\pi_5^2\pi_5'^2, 72.\\
\end{array}
\end{eqnarray*}

\renewcommand{\arraystretch}{1.0}

{\bf (7.3)} {\it If $[\sigma]$ denotes a conjugacy class above, then the fix lines of the elements in $[\sigma]$ form one orbit under $G_n$.}\\
\bprf
The statement is clear when $\sigma$ has just one line of fix points. We prove that the two fix lines of the elements in $[\sigma_{24}]$, $[\pi_3\pi_4\pi_3'\pi_4']$ or $[\pi_3\pi_4\sigma_4]$ are equivalent under $G_n$. The latter are eigenspaces of $\ci^4$ with eigenvalues $1$ and $-1$. Remember that if $\pi$ is an element in one of the previous classes then $-\pi$ is in the same conjugacy class too. It has the same eigenspaces but with eigenvalues interchanged. So we can find a matrix in $G_n$ which maps one line to the other and vice versa.
\eprf

\noindent
This shows in particular that every $G_n$-invariant property, which holds for a special fix line of an element in a conjugacy class above, holds for each other fix line of the elements in the same conjugacy class.\\  
{\bf (7.4)} We give now the generators of the fix line(s) of the representatives above, which we will need later, $\tau=\frac{1}{2}(1+\sqrt{5})$:
\begin{eqnarray*}
\begin{array}{l}
\sigma_{24}: <(0:0:1:0),(1:0:0:0)>, <(0:0:0:1),(0:1:0:0)>;\\ \pi_3\pi_3':<(1:0:0:0), (0:1:-1:1)>; \pi_3\pi_3'^2:<(0:1:1:0),(0:-1:0:1)>; \\\pi_3\pi_4\pi_3'\pi_4:<(1:0:0:0),(0:1:0:1)>, <(0:1:0:-1), (0:0:1:0)>; \\\pi_3\pi_4\sigma_4:<(1:\sqrt{2}:1:0), (0:1:\sqrt{2}:1)>, <(\sqrt{2}:-1:0:1), (1:-\sqrt{2}:1:0)>; \\\pi_5\pi_5':<(1:0:0:0), (0:0:\tau-1:1)>.
\end{array}
 \end{eqnarray*}
We can formulate a result about the intersection points of the fix lines:
{\bf (7.5)}{\it the intersection points of the previous lines are real. In particular they are not on the quadric.}\\
\bprf
Let $L$$\not=$$L'$  be some fix lines and let $x$ $\in L\cap L'$. Since the matrices of $G_n$ which fix the lines are real, it follows that $\bar x$$\in$$L\cap L'$ too. Hence $x=\bar x$.
\eprf

\paragraph{8. Coverings of $\piu$.}\label{coverings}
A pencil in $\pitr$ defines a morphism (away from the base locus)
\begin{eqnarray*}
\begin{array}{lll}
\pitr&\longrightarrow &\mathbb{P}_1\\
x&\mapsto&(S_n(x):Q_n(x)).
\end{array}
\end{eqnarray*}
Let $L$ be a fix line which does not meet the base locus. This morphism restricts on $L$ to a $n:1$ cover of $\mathbb{P}_1$, we call it $f$. If $L$ meets the base locus of $F_n(\lambda)$, then $f$ is not defined at the points of intersection $z_{13}$, $z_{24}$ (notation of paragraph {\bf 7.}). By {\bf (7.1)} the line $L$ meets each surface in $F_n(\lambda)$ with multiplicity one at these points. So $f$ extends to a cover 
\begin{eqnarray*}
\bar f:L\longrightarrow \mathbb{P}_1
\end{eqnarray*} 
of degree $n-2$, having branch points of order $\frac{n}{2}-1$ at $z_{13}$, $z_{24}$. Using Hurwitz's formula for the degree of the ramification locus of a curve morphism, we get\\
{\bf (8.1)} {\it the degree of the ramification locus of the morphism $f$ is $2n-2$ and of $\bar f$ it is $2n-6$ $(n=6,8,12)$.}\\
Of course the singular points of the surfaces in the pencil are ramification points of the previous morphisms. So the degree of the ramification locus gives an estimate for the number of singular points on the fix lines. In fact if we do not consider the points on the multiple quadric, we find that these points are at most $n$ if $L$ does not intersect the base locus, these are $n-2$ otherwise. 
\paragraph{9. Singular surfaces.}
The aim of this paragraph is to find  the singular surfaces in the pencil and their number of singular points.\\
{\bf (9.1)} {\it The pencils $F_6(\lambda)$ and $F_{12}(\lambda)$.} We have seen that the surfaces are invariant under the reflection groups of the 24-cell and of the 600-cell. We consider these in position as in \cite{cox} p. 157, where by the 600-cell and its reciprocal we interchange the $x_0$ and $x_1$ coordinate. The vertices of $\{3,4,3\}$ are the permutations of $(\pm 1, \pm 1, 0,0)$, those of $\{3,3,5\}$ are the permutations of $(\pm 2, 0,0,0)$, $(\pm 1,\pm 1, \pm 1,\pm 1)$ and the even permutations of $(\pm 1, \pm \tau, \pm \tau^{-1}, 0)$. Of course the $N_0$-vertices, $N_1$-edges, $N_2$-faces and $N_3$-cells form one orbit each under the reflection group action. We consider besides the vertices and the middle points of the edges, the vertices and the middle points of the edges of the reciprocal of $\{3,4,3\}$, which we denote by $\{3,4,3\}'$ and of the reciprocal of $\{5,3,3\}$ (their number is equal respectively to $N_3$ and $N_2$). We recall the coordinates of the vertices, which are given in \cite{cox} p. 156-157. For $\{3,4,3\}'$ these are the permutations of $(\pm 2,0,0,0)$ and $(\pm 1,\pm 1,\pm 1,\pm 1)$, for $\{5,3,3\}$ we have the permutations of $(\pm 2, \pm 2,0,0)$, $(\pm\sqrt{5},\pm1,\pm 1, \pm 1)$, $(\pm \tau, \pm \tau, \pm \tau, \pm \tau^{-2})$, $(\pm \tau^2,\pm \tau^{-1}, \pm \tau^{-1}, \pm \tau^{-1})$ and the even permutations of $(\pm \tau^{-2}, \pm \tau^{2}, \pm 1, 0)$, $(\pm\tau^{-1}, \sqrt{5}, \tau,0)$, $(\pm 1,\pm 2,\pm \tau, \pm \tau^{-1})$.  In these way we get four distinct orbits of points. If we consider these in $\mathbb{P}_3$, the orbits have half length, and we can multiply by some non zero scalar factor without changing the coordinates of the points. Substituting one point of each orbit in the equations of the pencil we get the values of $\lambda$, s.t. the corresponding surface contains the whole orbit. Since the partial derivatives of the equation vanish, these points are in fact singular (one needs to consider just one point in each orbit and to do the calculations with MAPLE). We give below the length of the orbits in $\pitr$ and the values of $\lambda$ of the  surfaces $F_6(\lambda)$, resp. $F_{12}(\lambda)$, which contain the orbits.

\renewcommand{\arraystretch}{1.6}

\begin{center}

\begin{tabular}{ccccc}
\multicolumn{5}{c}{24-cell}\\
orbit:&12&48&48&12\\
\hline
$\lambda$:&-1&$-\frac{2}{3}$&$-\frac{7}{12}$&$-\frac{1}{4}$\\
\end{tabular}~~
\begin{tabular}{ccccc}
\multicolumn{5}{c}{600-cell}\\
orbit:&300&600&360&60\\
\hline
$\lambda$:&$-\frac{3}{32}$&$-\frac{22}{243}$&$-\frac{2}{25}$&0
\end{tabular}

\end{center}

\renewcommand{\arraystretch}{1.0}
We show now that we have no more singular points and singular surfaces in the pencils. To do this we consider the lines of fix points for elements of $G_6$, resp. $G_{12}$. In fact, for our aim (cf. 1., 2., and 3. below), it is enough to consider just one line in each $G_n$-orbit and its intersection points with the above given singular surfaces. By {\bf (7.3)} the other lines meet the same surfaces in the same number of singular points and these form one $G_n$-orbit. We take the lines of {\bf (7.4)} and the fix line $<(0:0:1:0), (1:0:0:0)>$ of $\sigma_{24}$: 

\vspace*{0.3cm}

\renewcommand{\arraystretch}{1.0}

\begin{eqnarray*}
\begin{array}{c|c|c|c|c}
{\rm fix~line~of}&\multicolumn{2}{c|}{{\rm surface}}&\multicolumn{2}{c}{{\rm int.~points}}\\
\hline
\sigma_{24}&F_6(-1)&F_6(-\frac{1}{4})&(0:0:1:0),&(1:0:\pm 1:0)\\
&& &(1:0:0:0)&\\
\pi_3\pi_3'&F_6(-1)&F_6(-\frac{2}{3})&(\pm 1:1:-1:1), &(\pm 3:1:-1:1), \\
&&&(1:0:0:0)&(0:1:-1:1)\\
\pi_3\pi_3'^2&F_6(-\frac{7}{12})&F_6(-\frac{1}{4})&(0:2:1:-1),&(0:1:1:0),\\
&&& (0:-1:1:2),&(0:-1:0:1),\\
&&&(0:1:2:1)&(0:0:1:1)\\
\end{array}
\end{eqnarray*}

\renewcommand{\arraystretch}{1.0}


\vspace*{0.3cm}

In the case of the pencil $F_{12}(\lambda)$, we take the representative $\pi_5^2\sigma_2\pi_5'^2\sigma_4$  instead of $\pi_3\pi_3'$ (it is more convenient for the computations that we are going to do). Its fix line is $<(1:0:0:0),(0:\tau^{2}:1:0)>$:

\vspace*{0.3cm}

\begin{eqnarray*}
\begin{array}{c|c|c|c|c}
{\rm fix~line~of}&\multicolumn{2}{c|}{{\rm surface}}&\multicolumn{2}{c}{{\rm int.~points}}\\
\hline
\sigma_{24}&F_{12}(-\frac{3}{32})&F_{12}(-\frac{22}{243})&(1:0:\pm 1:0)&(1:0:\pm \tau^{2}:0),\\
&&&&(\pm \tau^{2}:0:1:0)\\
&F_{12}(-\frac{2}{25})&F_{12}(0)&(1:0:\pm\tau:0),&(1:0:0:0),\\
&&&(\pm \tau:0:1:0)&(0:0:1:0)\\
\pi_5^2\sigma_2\pi_5'^2\sigma_4&F_{12}(-\frac{3}{32})&F_{12}(-\frac{22}{243})&(\pm \tau^2:1:\tau^{-2}:0),&(0:\tau^2:1:0),\\
&&&(\pm \tau^{-2}:\tau^{2}:1:0),&(\pm 3:\tau:\tau^{-1}:0)\\
&&&(\pm\sqrt{5}:\tau:\tau^{-1}:0)&\\
&F_{12}(0)&&(1:0:0:0),&\\
&&&(\pm 1:\tau:\tau^{-1}:0)&\\
\pi_5\pi_5'&F_{12}(-\frac{2}{25})&F_{12}(0)&(0:0:1:\tau),&(1:0:0:0),\\
&&&(\pm\tau^2:0:\tau^{-1}:1),&(\pm\tau^{-1}:0:1:\tau),\\
&&&(\pm \sqrt{5}:0:\tau:\tau^2)&(\pm\tau:0:\tau^{-1}:1)\\
\end{array}
\end{eqnarray*}






\vspace*{0.3cm}

By the observation above and the tables follows that:
\begin{enumerate}
\item the points of intersection of the fix lines with the surfaces $F_6(-\frac{1}{4})$ and $F_6(-1)$, resp. $F_{12}(0)$ and $F_{12}(-\frac{3}{32})$  are in fact vertices of the 24-cell $\{3,4,3\}$, resp. of the reciprocal $\{3,4,3\}'$, and of the $\{3,3,5\}$ and its reciprocal $\{5,3,3\}$.
\item The points of intersection with the surfaces $F_6(-\frac{7}{12})$ and $F_6(-\frac{2}{3})$ are mid points of the edges of the $\{3,4,3\}$, resp. of the $\{3,4,3\}'$.
\item In the case of the two remaining surfaces $F_{12}(-\frac{2}{25})$ and $F_{12}(-\frac{22}{243})$ we have to do some more consideration. The middle point of the edge connecting the two neighboring vertices $(1,\tau,\pm\tau^{-1},0)$ is $(1,\tau,0,0)$. So all the permutations (with  + or - sign) of these points are middle points of the edges of the $\{3,3,5\}$. Similarly the middle points of the edge connecting the two points $(\pm \tau^{-2},\tau^{2},1,0)$ is $(0,\tau^{2},1,0)$, hence the permutations of these points (with + and - sign) are middle points of the edges of the $\{5,3,3\}$. So we find the intersection points of the fix line of $\sigma_{24}$ and the surfaces $F_{12}(-\frac{2}{25})$ and $F_{12}(-\frac{22}{243})$. The fix lines of the remaining two matrices $\pi_5^2\sigma_2\pi_5'^2\sigma_4$ and $\pi_5\pi_5'$ contain the points $(0:\tau^{2}:1:0)$$\in F_{12}(-\frac{22}{243})$ resp. $(0:0:1:\tau)$$\in F_{12}(-\frac{2}{25})$. The other points of intersection form one orbit under $\pi_5^2\sigma_2$, resp. $\pi_5$ which let the line invariant. Hence all the points of intersection are middle points of the edges of one polyhedron, resp. of the other.
\end{enumerate}
Now we are almost done. A calculation with MAPLE shows that the fix line of $\sigma_{24}$, resp. $\pi_5\pi_5'$ contains two points of the base locus of $F_6(\lambda)$, resp. of $F_{12}(\lambda)$. Using the estimations {\bf (8.1)} for the number of singular points on the fix lines we see that in fact we have found the maximal number possible. We conclude that
\begin{itemize}
\item we have no more singular surfaces and singular points in the pencils,
\item the points which we have found are all double points,
\item the singular points form one orbit already under the action of the groups $G_6$ and $G_{12}$. \\
\bprf
In the case of the surfaces $F_{12}(-\frac{2}{25})$ and $F_{12}(-\frac{22}{243})$, the assertion follows by 3. above and {\bf (7.3)}. In the case of the surfaces $F_{12}(-\frac{3}{32})$ and $F_{12}(0)$ we proceed in a similar way. We consider the two points of intersection of each surface with the fix line of $\sigma_{24}$ and we observe that the matrix $\sigma_2$ maps one point to the other. Hence we conclude again by {\bf (7.3)}.
\eprf
\end{itemize}
{\bf (9.2)} {\it The pencil $F_8(\lambda)$.} As we remarked in paragraph {\bf 3.}, the group $G_8$ maps the vertices of a $\{3,4,3\}$ in those of a $\{3,4,3\}'$, (the vertices are taken now with coordinates as in {\bf 3.}). The same holds for the middle points of the edges. In fact a computation with MAPLE shows that these two orbits belong respectively to the surfaces $F_8(-1)$ and $F_8(-\frac{5}{9})$ and the points are singular. Obviously these are the points of the surfaces $F_6(-1)$ and $F_6(-\frac{1}{4})$ together, resp. of $F_6(-\frac{2}{3})$ and $F_6(-\frac{7}{12})$ together. We look now for some more orbits of points. For instance consider the middle points of the segments connecting the vertices of $\{3,4,3\}$ with its reciprocal $\{3,4,3\}'$. These are the permutations of:
\begin{eqnarray*}
\begin{array}{l}
(1)~(\pm \frac{\sqrt{2}+1}{2},\pm\frac{1}{2},0,0),  (\pm \frac{\sqrt{2}-1}{2},\pm\frac{1}{2},0,0);\\ (2)~\frac{1}{2}(\pm (1+\frac{\sqrt{2}}{2}),\pm (1+\frac{\sqrt{2}}{2}),\pm\frac{\sqrt{2}}{2},\pm\frac{\sqrt{2}}{2}), \frac{1}{2}(\pm (1-\frac{\sqrt{2}}{2}),\pm (1-\frac{\sqrt{2}}{2}),\pm\frac{\sqrt{2}}{2},\pm\frac{\sqrt{2}}{2});\\ (3)~\frac{1}{2}(\pm 1,\pm 1,\pm \sqrt{2},0); \\(4)~\frac{1}{2}(\pm (1+\frac{\sqrt{2}}{2}),\pm (1-\frac{\sqrt{2}}{2}),\pm\frac{\sqrt{2}}{2},\pm\frac{\sqrt{2}}{2}).
\end{array}
\end{eqnarray*}
Again a computation with MAPLE shows that the points (1), (2), which now we consider in $\pitr$, are singular on the surface $F_8(-\frac{3}{4})$ and the points (3), (4), in $\pitr$, are singular on $F_8(-\frac{9}{16})$. We take, in fact just one point in each group. Then because of symmetry reasons the points in the same group are singular on the same surface. We want to understand in how many $G_8$-orbits do these points split. We know already that we have at least two orbits since the points  belong to two different surfaces. To proceed we use the lines of fix points. As in {\bf (9.1)} we give the intersections of the fix lines of {\bf (7.4)} with the found singular surfaces. For $\sigma_{24}$ we take the same fix line as in {\bf (9.1)}, for $\pi_3\pi_4\sigma_4$ we take $<(1:\sqrt{2}:1:0),(0:1:\sqrt{2}:1)>$ and for $\pi_3\pi_4\pi_3'\pi_4'$ we take $<(1:0:0:0), (0:1:0:1)>$. Put $a:=1+\sqrt{2}$.

\vspace*{0.3cm}

\begin{eqnarray*}
\begin{array}{c|c|c|c|c}
{\rm fix~line~of}&\multicolumn{2}{c|}{{\rm surface}}&\multicolumn{2}{c}{{\rm int.~points}}\\
\hline
\sigma_{24}&F_8(-1)&F_8(-\frac{3}{4})&(0:0:1:0),&(1:0:\pm a),\\
&&& (1:0:0:0),&(\pm a:0:1:0)\\
&&&(1:0:\pm 1:0)&\\
\pi_3\pi_3'&F_8(-1)&F_8(-\frac{5}{9})&(\pm 1:1:-1:1), &(\pm 3:1:-1:1),\\
&&&(1:0:0:0)& (0:1:-1:1)\\
\pi_3\pi_4\sigma_4&F_8(-\frac{3}{4})&F_8(-\frac{9}{16})&(a:1:-1:-a),&(1:\sqrt{2}:1:0),\\
&&&(1:a:a:1),&(0:1:\sqrt{2}:1),\\
&&&(-1:1:a:a),&(1:0:-1:-\sqrt{2}),\\
&&&(a:a:1:-1)&(\sqrt{2}:1:0:-1)\\
\pi_3\pi_4\pi_3'\pi_4'&F_8(-1)&F_8(-\frac{5}{9})&(1:0:0:0),&(\pm 1:1:0:1),\\
&&&(0:1:0:1)& (\pm 2:1:0:1)\\
&F_8(-\frac{9}{16})&&(\pm \sqrt{2}:1:0:1)&\\
\end{array}
\end{eqnarray*}




\vspace*{0.3cm}

We show that the points (1), (2) on $F_8(-\frac{3}{4})$, resp. the points (3),  (4) on $F_8(-\frac{9}{16})$ form one $G_8$-orbit.\\
\bprf
By the tables above the fix line of $\sigma_{24}$ and of $\pi_3\pi_4\sigma_4$ meets the surface $F_8(-\frac{3}{4})$ at points which are in (1), resp.  (2) (we have just to multiply the points by some scalar factor). Hence we cannot have other singular points besides the points (1), (2) together. We have to show:\\
(a) the points on the fix lines form one orbit under the fix group of the line,\\
(b) the fix lines of the elements in $[\sigma_{24}]$ meet some fix line of the elements in $[\pi_3\pi_4\sigma_4]$ at the singular points of $F_8(-\frac{3}{4})$.\\
We remarked in paragraph {\bf 7.}, {\bf (7.2)}, that a matrix of order two in $[\sigma_2]$ is the square of a matrix of order four in $[\pi_4]$. This holds clearly for $\sigma_4$ too. Put $\sigma_4=\gamma^2$, with $\gamma\in[\pi_4]$, then $\gamma$ and $\sigma_4$ commute, hence the four points on the fix lines of $\sigma_{24}$, resp. $\pi_3\pi_4\sigma_4$ form one orbit under $\gamma$. This shows (a). Consider now the fix line $<(1:0:0:1), (0:1:1:0)>$ of the element $\sigma_2\sigma_3$$\in$$[\sigma_{24}]$, it contains the point $(1:a:a:1)$ of the fix line of $\pi_3\pi_4\sigma_4$. This shows (b) and completes the case of $F_8(-\frac{3}{4})$. We do the same for the points (3), (4) on $F_8(-\frac{9}{16})$. We consider the fix line of $\pi_3\pi_4\sigma_4$ and of $\pi_3\pi_4\pi_3'\pi_4'$. These meet $F_8(-\frac{9}{16})$ in four resp. two points, an argumentation as before shows that they form one orbit under the fix group of the line, so we get (a). To show (b), we consider the fix line $<(0:0:0:1), (-1:0:1:0)>$ of $\sigma_2\pi_3\pi_4\pi_3'\pi_4'$$\in$$[\pi_3\pi_4\pi_3'\pi_4']$ it contains the point $(1:0:-1:-\sqrt{2})$ of the fix line of $\pi_3\pi_4\sigma_4$ and we have done.
\eprf 
A calculation with MAPLE shows that the fix line of $\pi_3\pi_3'$ meets the base locus of $F_8(\lambda)$ in two points. Using the estimations for the number of singular points on the fix lines we see that:
\begin{itemize}
\item we have no more singular surfaces and singular points in the pencil $F_8(\lambda)$,
\item the singular points are all double points,
\item the singular surfaces have just one orbit of singular points each. We resume below their length in $\pitr$ and the values of $\lambda$ s.t. $F_8(\lambda)$ is singular.

\renewcommand{\arraystretch}{1.6}

\begin{eqnarray*}
\begin{array}{ccccc}
{\rm orbit}:&24&72&144&96\\
\hline
\lambda:&-1&-\frac{3}{4}&-\frac{9}{16}&-\frac{5}{9}\\
\end{array}
\end{eqnarray*}

\renewcommand{\arraystretch}{1.0}

\end{itemize} 
\paragraph{10. The singular points are nodes.} We take now an arbitrary singular point in each orbit and we compute the Hessian matrix at this point. A calculation with MAPLE shows that in any case the matrix has rank three, hence the point is a node and so are all the points in its orbit.

\paragraph{11. Configurations of lines and points.} We recall the following\\
{\bf (11.1) Definition.} A {\it space configuration} of lines and points is a system of $l$ lines and $p$ points s.t. each line contains $\pi$ of the given points and each point belongs to $\lambda$ lines. We say that we have a $(p_{\lambda},l_{\pi})$ configuration. Moreover we have the formula
\begin{eqnarray*}
\begin{array}{c}
p\cdot\lambda=l\cdot\pi.
\end{array}
\end{eqnarray*}
The fix lines of elements in the same conjugacy class and the singular points on a singular surface form a space configuration of lines and points.
By the formula above and the results of paragraph {\bf 10.}, we can compute how many fix lines of the same orbit pass through each singular point. We list the configurations below. In the first column we give a representative of the conjugacy class which we consider.

\vspace*{0.3cm}

\renewcommand{\arraystretch}{1.3}

\begin{eqnarray*}
\begin{array}{c|c|c|c|c}
{\rm repr.~of~the}&\multicolumn{2}{c|}{{\rm sing.~points~of }}&\multicolumn{2}{c}{{\rm configuration}}\\
{\rm conj.~class}&\multicolumn{2}{c|}{{\rm }}&\multicolumn{2}{c}{}\\
\hline
\sigma_{24}&F_{6}(-1)&F_{6}(-\frac{1}{4})&(18_2,12_3)&(18_2,12_3)\\
&F_{8}(-1)&F_{8}(-\frac{3}{4})&(18_4,24_3)&(18_4,72_1)\\
&F_{12}(-\frac{3}{32})&F_{12}(-\frac{22}{243})&(450_2,300_3)&(450_4,600_3)\\
&F_{12}(-\frac{2}{25})&F_{12}(0)&(450_4,360_5)&(450_2,60_{15})\\
\pi_3\pi_3'&F_{6}(-1)&F_{6}(-\frac{2}{3})&(16_3,12_4)&(16_3,48_1)\\
&F_{8}(-1)&F_{8}(-\frac{5}{9})&(32_3,24_4)&(32_3,96_1)\\
&F_{12}(-\frac{3}{32})&F_{12}(-\frac{22}{243})&(200_6,300_4)&(200_3,600_1)\\
&F_{12}(0)&&(200_3,60_{10})&\\
\pi_3\pi_3'^2&F_{6}(-\frac{7}{12})&F_{6}(-\frac{1}{4})&(16_3,48_1)&(16_3,12_4)\\
\pi_3\pi_4\sigma_4&F_{8}(-\frac{3}{4})&F_{8}(-\frac{9}{16})&(36_4,72_2)&(36_4,144_1)\\

\pi_3\pi_4\pi_3'\pi_4'&F_{8}(-1)&F_{8}(-\frac{9}{16})&(72_2,24_6)&(72_2,144_1)\\
&F_{8}(-\frac{5}{9})&&(72_4,96_3)&\\
\pi_5\pi_5'&F_{12}(-\frac{2}{25})&F_{12}(0)&(72_5,360_1)&(72_5,60_6)\\
\end{array}
\end{eqnarray*}



\renewcommand{\arraystretch}{1.0}

\noindent
It is interesting to observe that the singular points are intersection points of the fix lines, hence {\bf (7.5)} confirms that they are all real points.

\newpage

\paragraph{12. Computer picture.}
We exhibit a computer picture of the surface of degree 12 with 600 nodes. This has been realized by the program SURF written by S. Endra\ss.
%
%
%
\begin{center}
\hspace*{-4cm}
\epsfig{file=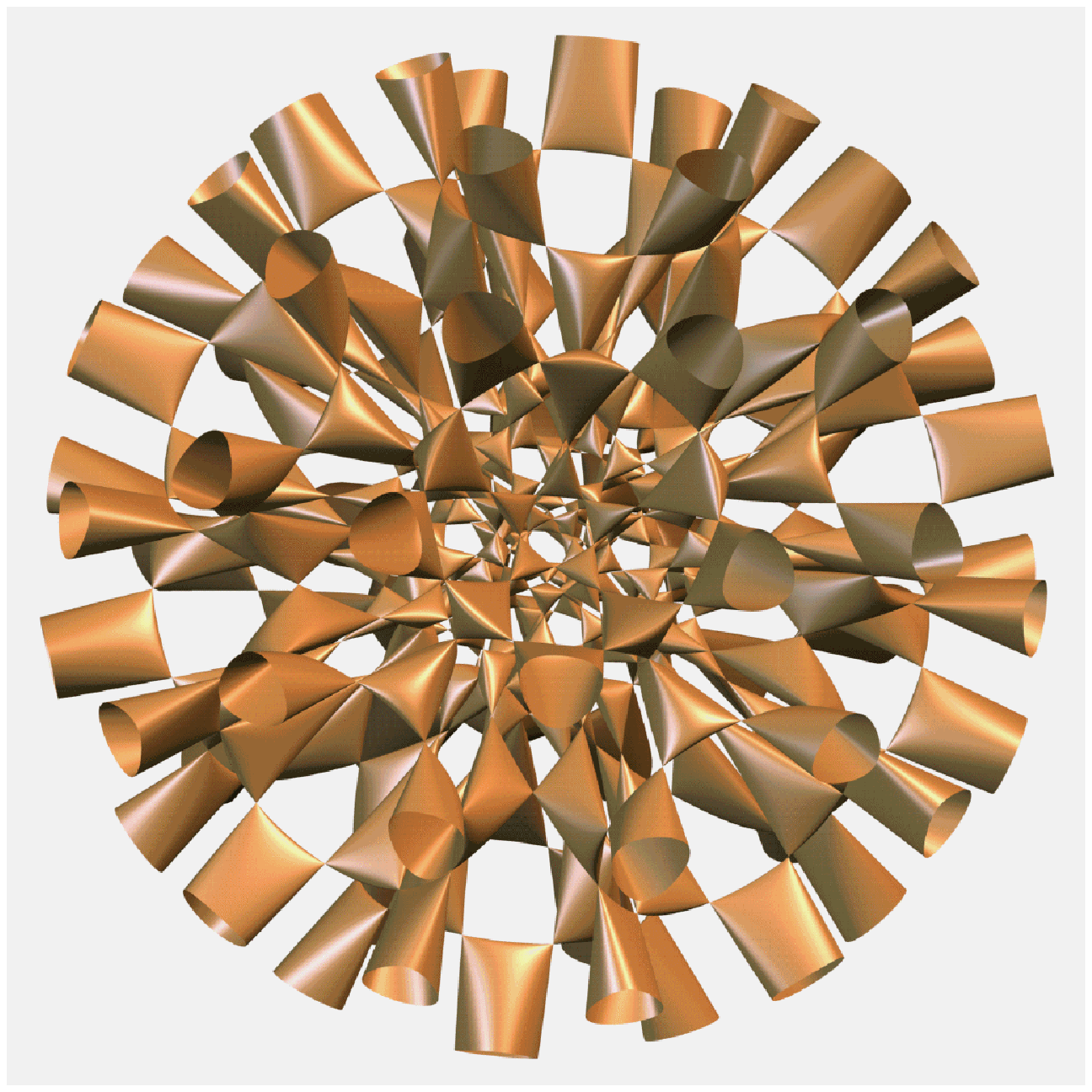}
\hspace*{-4cm}\ \\[3ex]

{$A_5\times A_5$--symmetric surface of degree 12}\\[0.5cm]
{with 600 nodes}

\end{center}
%


\cleardoublepage


\begin{thebibliography}{99}
\bibitem{Ben}
Benson, D. J.: {\it Polynomial Invariants of Finite Groups}, London Math. Society LNS 190, Cambridge University Press (1993).
\bibitem{burniside}
Burnside, W.:{\it Theory of groups of finite order}, Dover Publications, Inc. (1955).
\bibitem{cox0}
Coxeter, H. S. M.:{\it The product of the generators of a finite group generated by reflections}, Duke Math. J. Vol. 18 (1951) 765-782.
\bibitem{cox}
Coxeter, H. M. S.: {\it Regular polytopes (second edition)}, The Macmillan company, New York (1963).
\bibitem{hudson}
Hudson, R. W. H. T.: {\it Kummer's Quartic Surface}, Cambridge University Press (1905).
\bibitem{jessop}
Jessop, C. M.: {\it A Treatise on the Line Complex}, Chelsea Publishing Company, New York (1969).
\bibitem{klein}
Klein, F.: {\it Vorlesungen \"uber das Ikosaeder und die Aufl\"osung der Glei\-chungen vom f\"unften Grade}, Nachdr. der Ausg. Leipzig, Teubner 1884, hrsg. mit einer Einf\"uhrung und mit Kommentaren von Peter Slodowy, Birkh\"auser-B. G. Teubner (1993).
\bibitem{kreiss}
Kreiss, H.O.: {\it \"Uber syzygetische Fl\"achen}, Ann. di Matematica (2) 41 (1955) 105-111.
\bibitem{miyaoka}
Miyaoka, Y.: {\it The maximal Number of Quotient Singularities on Surfaces with Given Numerical Invariants}, Math. Ann. 268 (1984) 159-171. 
\bibitem{mukai}
Mukai, S.: {\it Moduli of Abelian  surfaces, and regular polyhedral groups}, Proceedings of symposium ``Moduli of algebraic varieties'' Sapporo, January 1999.
\bibitem{racah}
Racah, G.: {\it Sulla caratterizzazione delle rappresentazioni irriducibili dei gruppi semisemplici di Lie}, Rend. Acad. Naz. dei Lincei, Classe di Scienze fisiche, matematiche e naturali (8), vol. 8 (1950) 108-112.
\bibitem{tits}
Tits, J.: {\it Liesche Gruppen und Algebren}, Springer (1983).
\bibitem{var}
Varchenko, A. N.: {\it On Semicontinuity of the Spectrum and an Upper Estimate for the Number of Singular Points of a Projective Hypersurface}, Soviet. Math. Dokl. Vol. 27 (1983) No. 3 735-739.  
\end{thebibliography}
\end{document}